\documentclass[12pt,thmsa]{article}%
\usepackage{amsfonts}
\usepackage{sw20bams}
\usepackage{amsmath}
\usepackage{amssymb}
\usepackage{graphicx}%
\providecommand{\U}[1]{\protect\rule{.1in}{.1in}}
\begin{document}

\author{Steven R. Finch}
\title{Covering a Sphere with Four Random Circular Caps}
\date{February 24, 2015}
\maketitle

\begin{abstract}
Let $p(\omega)$ denote the probability that four random circular caps of
angular radius $70^{\circ}<\omega<90^{\circ}$ cover the unit sphere $S^{2} $.
\ An exact expression for $p(\omega)$ is unknown. \ We give nontrivial lower
bounds for $p(\omega)$ when $\omega>84^{\circ}$; no improvement on the
inequality $p(\omega)\geq0$ for $\omega<84^{\circ}$ is yet feasible. \ A\ dual
problem involving randomly inscribed well-centered tetrahedra is also examined.

\end{abstract}

\footnotetext{Copyright \copyright \ 2015 by Steven R. Finch. All rights
reserved.}Let $S^{2}$ denote the two-dimensional sphere in $\mathbb{R}^{3}$ of
unit radius. \ Let $X_{1}$, $X_{2}$, $X_{3}$, $X_{4}$ be four points that are
independent and uniformly distributed over the surface of $S^{2}$. Each
$X_{j}$ is taken to be the center of a circular cap of angular radius $\omega
$:%
\[%
\begin{array}
[c]{ccc}%
\left\{  Y\in S^{2}:\arccos(X_{j}\cdot Y)\leq\omega\right\}  , &  &
0\leq\omega\leq\pi
\end{array}
\]
where $X\cdot Y$ is the usual inner product between vectors in $\mathbb{R}%
^{3}$. \ A point, hemisphere and $S^{2}$ itself are the outcomes if $\omega
=0$, $\pi/2$ and $\pi$. The sphere is said to be \textbf{covered} if each
point of $S^{2}$ belongs to at least one cap. It is known that the coverage
probability $p(\omega)$ satisfies \cite{Mi2, BCL}%
\[
p(\omega)=\left\{
\begin{array}
[c]{lll}%
0 &  & \text{if }0\leq\omega\leq\omega_{0},\\
1-2\cos^{6}\left(  \frac{\omega}{2}\right)  \left[  8-9\cos^{2}\left(
\frac{\omega}{2}\right)  \right]   &  & \text{if }\frac{\pi}{2}\leq\omega
\leq\pi
\end{array}
\right.
\]
\ where%
\[
\omega_{0}=\arccos\left(  \tfrac{1}{3}\right)  =1.23095941...\approx
70.53^{\circ}.
\]
Finding an exact expression for $p(\omega)$ in the interval $\omega_{0}%
<\omega<\pi/2$ remains an open problem. \ Early works in this area include
\cite{MFG, Gbt}; computer simulations appear in \cite{SoS, Zhg} and historical
overviews in \cite{Smn, Hal, BCL}.\ All such papers deal with $n$ circular
caps, where $n$ is often large. \ Our restriction to the case $n=4$ is
somewhat unusual, as far as is known.

Gilbert \cite{Gbt} gave the following bounds on $p(\omega)$:
\[
0\leq p(\omega)\leq%
{\displaystyle\sum\limits_{k=0}^{2}}
\tbinom{2}{k}(-1)^{k}\left(  1-k\sin^{2}\left(  \tfrac{\omega}{2}\right)
\right)  ^{4}=1-2\cos^{8}\left(  \tfrac{\omega}{2}\right)  +\cos^{4}(\omega)
\]
(his left-hand side is negative for $n=4$ and thus replaced by zero). The
upper bound is weak but not relevant to us. We will focus on the lower bound.
Our improvement, based on elaborate formulation, occurs only when
$\omega>84.25^{\circ}$. \ For example, if $\omega=88^{\circ}$, then
$0.0765<p(\omega)<0.8567$ via our work and $p(\omega)\approx0.08$ via computer
simulation. \ 

Dual to the coverage problem is the following. \ Given $A$, $B$, $C$, $D$ on
$S^{2}$, there is a unique circular cap of minimal angular radius
$\theta_{\min}$ containing all four points. \ If $A$, $B$, $C$, $D$ are
independent and uniform on $S^{2}$, then the cumulative distribution function
$\Phi$ for the random variable $\theta_{\min}$ satisfies \cite{Mi2, BCL}%
\[%
\begin{array}
[c]{ccc}%
1-\Phi(\pi-\omega)=p(\omega), &  & 0\leq\omega\leq\pi.
\end{array}
\]
Define a distribution%
\[%
\begin{array}
[c]{ccc}%
F(\xi)=8\left(  \Phi(\xi)-\tfrac{7}{8}\right)  =8\,\mathbb{P}\left\{
\frac{\pi}{2}<\theta_{\min}\leq\xi\right\}  , &  & \frac{\pi}{2}<\xi<\pi.
\end{array}
\]
Calculating $F(\xi)$ is hopeless; therefore we study a different variable
$\theta_{\operatorname*{abc}}$ with distribution $G$ that is computationally
feasible and satisfies $c\,G(\xi)\geq F(\xi)$ for some constant $c>0$.

Any three points from $\{A,B,C,D\}$ almost surely form both a chordal triangle
(with sides as straight lines through the interior of $S^{2}$) and a spherical
triangle (with sides as great circles on the surface of $S^{2}$). \ 

The chordal triangle $ABC$ determines a plane that partitions $S^{2}$ into two
complementary caps. Let $ABCD$ denote the tetrahedron inscribed in $S^{2} $
with base $ABC$ and apex $D$. Such a tetrahedron is called
\textbf{well-centered} if it contains the origin \cite{VrZ}. \ 

The hypothesis $\pi/2<\theta_{\min}$ implies that $ABCD$ is well-centered and,
without loss of generality, $ABC$ is acute. Further,$\ D$ is in the larger of
the two complementary caps and the angular radius $\theta_{\operatorname*{abc}%
}$ of the cap satisfies $\pi/2<\theta_{\operatorname*{abc}}$. \ The converse
is not true. Our main objective is to find the conditional probability density
function $g$ for $\theta_{\operatorname*{abc}}$ and the associated
distribution $G$. Since $ABC$ is one of four possible choices of base, we have
a crude estimate $4G(\xi)\geq F(\xi)$. \ 

While finding $g$, we compute the probability%
\[
\kappa=\mathbb{P}\left\{  ABCD\text{ is well-centered \& }ABC\text{ is
acute}\right\}  =0.10191818...
\]
which extends a well-known result $\mathbb{P}\left\{  ABCD\text{ is
well-centered}\right\}  =1/8$ from \cite{WiW, Wen, KAL, HwS, Ber} and another
result $\mathbb{P}\left\{  ABC\text{ is acute}\right\}  =1/2$ from \cite{Mi2}.
\ Details underlying the acuteness probability are given in the next section.
\ The expected number of acute faces of $ABCD$, given well-centeredness, is
$32\kappa$; a refined estimate is hence $(3.26138193...)G(\xi)\geq F(\xi)$. \ 

\section{Random Triangles of Circumradius $r$}

For comparison's sake, consider triangles formed by selecting three
independent points on the circle $r\,S^{1}$, where the radius $r$ is fixed.
The bivariate density for two arbitrary angles $\alpha$, $\beta$ is
\cite{Fin}
\[
\left\{
\begin{array}
[c]{lll}%
\dfrac{2}{\pi^{2}} &  & \text{if }0<\alpha<\pi\text{, }0<\beta<\pi\text{ and
}\alpha+\beta<\pi,\\
0 &  & \text{otherwise}%
\end{array}
\right.
\]
thus the acuteness probability for such triangles is $1/4$. \ The bivariate
density for two arbitrary sides $a$, $b$ is
\[
\left\{
\begin{array}
[c]{lll}%
\dfrac{4}{\pi^{2}r^{2}}\dfrac{1}{\sqrt{4-\left(  \frac{a}{r}\right)  ^{2}%
}\sqrt{4-\left(  \frac{b}{r}\right)  ^{2}}} &  & \text{if }0<\frac{a}%
{r}<2\text{ and }0<\frac{b}{r}<2,\\
0 &  & \text{otherwise;}%
\end{array}
\right.
\]
the sides $a$, $b$ are independent despite the fact that angles $\alpha$,
$\beta$ are \textit{dependent} and $a=2r\sin(\alpha)$, $b=2r\sin(\beta)$.

Consider now chordal triangles formed by selecting three independent points on
the sphere $S^{2}$. Given triangle $ABC$, let $r$ denote the radius of the
unique circle passing through vertices $A$, $B$, $C$. \ Unlike before, $r $ is
a random variable. Intuition might suggest that the density of $\alpha$,
$\beta$ and of $a/r$, $b/r$ should be similar to the preceding. \ In fact,
however, the bivariate density of $\alpha$, $\beta$ is%
\[
\left\{
\begin{array}
[c]{lll}%
\dfrac{8}{3\pi}\sin(\alpha)\sin(\beta)\sin(\alpha+\beta) &  & \text{if
}0<\alpha<\pi\text{, }0<\beta<\pi\text{ and }\alpha+\beta<\pi,\\
0 &  & \text{otherwise}%
\end{array}
\right.
\]
thus the acuteness probability for such triangles is $1/2$. \ The trivariate
density of $a$, $b$, $r$ is%
\[
\frac{a\,b}{6\pi r^{4}}\left\{  \frac{\frac{a}{r}\sqrt{4-\left(  \frac{b}%
{r}\right)  ^{2}}+\frac{b}{r}\sqrt{4-\left(  \frac{a}{r}\right)  ^{2}%
}+\left\vert \frac{a}{r}\sqrt{4-\left(  \frac{b}{r}\right)  ^{2}}-\frac{b}%
{r}\sqrt{4-\left(  \frac{a}{r}\right)  ^{2}}\right\vert }{\sqrt{4-\left(
\frac{a}{r}\right)  ^{2}}\sqrt{4-\left(  \frac{b}{r}\right)  ^{2}}}\right\}
\]
for $0<a<2r$ and $0<b<2r$. \ This is a $(2/3,1/3)$-weighted mixture of
densities. \ The portion coming from acute chordal triangles is%
\[
\frac{a\,b}{6\pi r^{4}}\frac{\frac{a}{r}\sqrt{4-\left(  \frac{b}{r}\right)
^{2}}+\frac{b}{r}\sqrt{4-\left(  \frac{a}{r}\right)  ^{2}}}{\sqrt{4-\left(
\frac{a}{r}\right)  ^{2}}\sqrt{4-\left(  \frac{b}{r}\right)  ^{2}}}%
=\frac{a\,b}{6\pi r^{4}}\left\{  \frac{a}{\sqrt{4r^{2}-a^{2}}}+\frac{b}%
{\sqrt{4r^{2}-b^{2}}}\right\}
\]
for $0<a<2r$, $0<b<2r$ and $a^{2}+b^{2}>4r^{2}$. \ We call this portion
$\delta(a,b,\theta)$ for future reference, where $r=\sin(\theta)$.

Proof of the above follows via a Jacobian determinant calculation based on
formula (4.18) of \cite{Mi1}; see also Theorem 3.2 of \cite{Mi2}. This is one
of several crucial contributions by Miles to our study.

\section{Conditional Probabilities and Angular Radius}

Let $\mathcal{E}$ denote the event $\{ABCD$ is well-centered \& $ABC$ is
acute$\}$. \ Following the proof of Theorem 3.3 of \cite{Mi2}, we have%
\[
\mathbb{P}\{\theta<\theta_{\operatorname*{abc}}<\theta+d\theta\}\cdot
\mathbb{P}\{\mathcal{E\;}|\;\theta\}=\mathbb{P}\{\theta<\theta
_{\operatorname*{abc}}<\theta+d\theta\mathcal{\;}|\;\mathcal{E}\}\cdot
\mathbb{P}\{\mathcal{E}\}
\]
where $\pi/2<\theta<\pi$. \ Formula (3.3) of \cite{Mi2} yields%
\[
\mathbb{P}\{\theta<\theta_{\operatorname*{abc}}<\theta+d\theta\}=\tfrac{3}%
{2}\sin^{3}(\theta)\,d\theta
\]
and finding
\[
\mathbb{P}\{\theta<\theta_{\operatorname*{abc}}<\theta+d\theta\mathcal{\;}%
|\;\mathcal{E}\}=g(\theta)\,d\theta
\]
is our main objective. \ We turn attention to $\mathbb{P}\{\mathcal{E\;}%
|\;\theta\}$.

Fix a chordal triangle $ABC$ for consideration. \ Let $\Delta$ denote the
spherical triangle with vertices $A$, $B$, $C$. The set of all points $D$ for
which the tetrahedron $ABCD$ is well-centered is $\Delta^{\prime}$, the
spherical triangle with vertices $A^{\prime}$, $B^{\prime}$, $C^{\prime}$
antipodal to $A$, $B$, $C$. \ This pictured in Figure 2 of \cite{VrZ} and is
stated in \cite{Ber}. Note that the area of $\Delta^{\prime}$ is equal to the
area of $\Delta$. \ 

Letting $ABC$ vary over all acute chordal triangles of fixed angular radius
$\pi/2<\theta<\pi$, the random quantity $\operatorname*{area}(\Delta)$ can be
as large as $-6\operatorname{arccot}\left(  \sqrt{3}\cos\theta\right)  -\pi$.
\ This is the area of an equilateral spherical triangle with vertices on a
circle of radius $\sin(\theta)$ \cite{BoE}. \ A formula for the mean area, as
opposed to the maximum area, turns out to be quite complicated. We derive this
in the next section.

\section{Kolmogorov-Robbins Theorem}

The desired probability $\mathbb{P}\{\mathcal{E\;}|\;\theta\}$ is the expected
value of $\operatorname*{area}(\Delta)/(4\pi)$, where random spherical
triangles $\Delta$ correspond to acute chordal triangles $ABC$ of fixed
angular radius $\pi/2<\theta<\pi$ \cite{Kgv, Rbn, KdM, Mae}.

Let $r=\sin(\theta)$ for convenience. If angles $\alpha$, $\beta$ of an acute
chordal triangle are known, then the remaining angle $\gamma$ satisfies%
\begin{align*}
\cos(\gamma)  &  =\cos\left(  \pi-\arcsin\left(  \tfrac{a}{2r}\right)
-\arcsin\left(  \tfrac{b}{2r}\right)  \right) \\
&  =\frac{a\,b-\sqrt{4r^{2}-a^{2}}\sqrt{4r^{2}-b^{2}}}{4r^{2}}.
\end{align*}
The corresponding spherical triangle has angle $\widetilde{\gamma}$ given by
\cite{Pea}
\begin{align*}
\widetilde{\gamma}  &  =\arccos\left(  \frac{4\cos(\gamma)-a\,b}{\sqrt
{4-a^{2}}\sqrt{4-b^{2}}}\right) \\
&  =\arccos\left(  \frac{(1-r^{2})a\,b-\sqrt{4r^{2}-a^{2}}\sqrt{4r^{2}-b^{2}}%
}{r^{2}\sqrt{4-a^{2}}\sqrt{4-b^{2}}}\right)  .
\end{align*}
Call this expression $\lambda(a,b,\theta)$. \ While $\operatorname*{area}%
(\Delta)=\widetilde{\alpha}+\widetilde{\beta}+\widetilde{\gamma}-\pi$ and
analogous expressions for $\widetilde{\alpha}$, $\widetilde{\beta}$ are
possible, it is simpler to employ $\mathbb{E}(\operatorname*{area}%
(\Delta))=3\,\mathbb{E}(\widetilde{\gamma})-\pi$. This yields \
\[
\mathbb{P}\{\mathcal{E\;}|\;\theta\}=%
{\displaystyle\int\limits_{0}^{2r}}
{\displaystyle\int\limits_{\sqrt{4r^{2}-b^{2}}}^{2r}}
\dfrac{1}{4\pi}\left(  3\lambda(a,b,\theta)-\pi\right)  \delta(a,b,\theta
)\,da\,db
\]
and%
\[
\mathbb{P}\{\mathcal{E}\}=%
{\displaystyle\int\limits_{\pi/2}^{\pi}}
\,\dfrac{3}{2}\sin^{3}(\theta)\,\mathbb{P}\{\mathcal{E\;}|\;\theta
\}\,d\theta=0.10191818...=\kappa;
\]
therefore we conclude that%
\[
G(\theta)=%
{\displaystyle\int\limits_{\pi/2}^{\theta}}
g(t)\,dt=%
{\displaystyle\int\limits_{\pi/2}^{\theta}}
\,\dfrac{3}{2\kappa}\sin^{3}(t)\,\mathbb{P}\{\mathcal{E\;}|\;t\}\,dt.
\]
In Figure 1, a plot of the density $g(\theta)$ is superimposed on a histogram
of simulated $\theta_{\operatorname*{abc}}$-values, $\theta\in\lbrack\pi
/2,\pi]$, based on $10^{6}$ random well-centered tetrahedra $ABCD$ with acute
faces $ABC$. \ The fit is excellent. First and second moments are also
indicated. \ All required integrals were evaluated via numerical methods
because no simplification via computer algebra seems likely.

Of course, our real interest is in $\theta_{\min}$, not $\theta
_{\operatorname*{abc}}$. \ Define a function%
\[
\psi(\theta)=\left\{
\begin{array}
[c]{lll}%
24\sin^{5}\left(  \tfrac{\theta}{2}\right)  \cos\left(  \tfrac{\theta}%
{2}\right)  \left[  2-3\sin^{2}\left(  \tfrac{\theta}{2}\right)  \right]  &  &
\text{if }0\leq\theta\leq\tfrac{\pi}{2},\\
\tfrac{1}{8}(32\kappa)g(\theta) &  & \text{if }\tfrac{\pi}{2}<\theta\leq\pi.
\end{array}
\right.
\]
In Figure 2, a plot of $\psi(\theta)$ is superimposed on a histogram of
simulated $\theta_{\min}$-values, $\theta\in\lbrack0,\pi]$, based on $10^{6}$
random (arbitrary) tetrahedra $ABCD$. \ The curve $\psi$ is represented by the
solid line and the fit is excellent on $[0,\pi/2]$. \ The area under $\psi$
over $[0,\pi/2]$ is $7/8$; the area under $\psi$ over $[\pi/2,\pi]$ is
$4\kappa=0.40767274...$. It visibly dominates the (unknown) density function
$\varphi$ for $\theta_{\min}$. \
\begin{figure}[ptb]%
\centering
\includegraphics[
height=5.3636in,
width=5.3636in
]%
{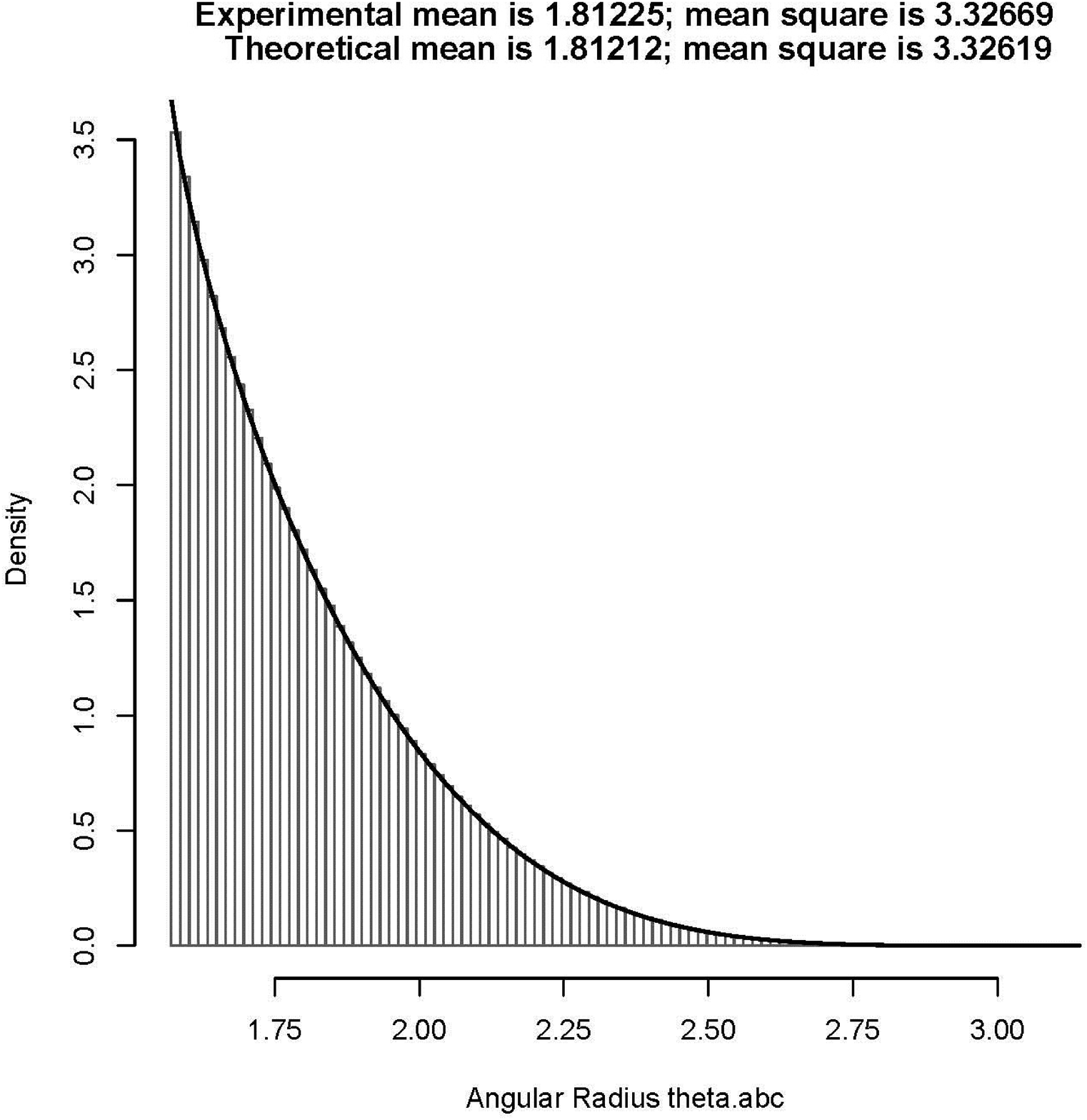}%
\caption{Experimental histogram of $\theta_{\operatorname*{abc}}$-values and
theoretical fit.}%
\end{figure}
%

\begin{figure}[ptb]%
\centering
\includegraphics[
height=5.3117in,
width=5.3117in
]%
{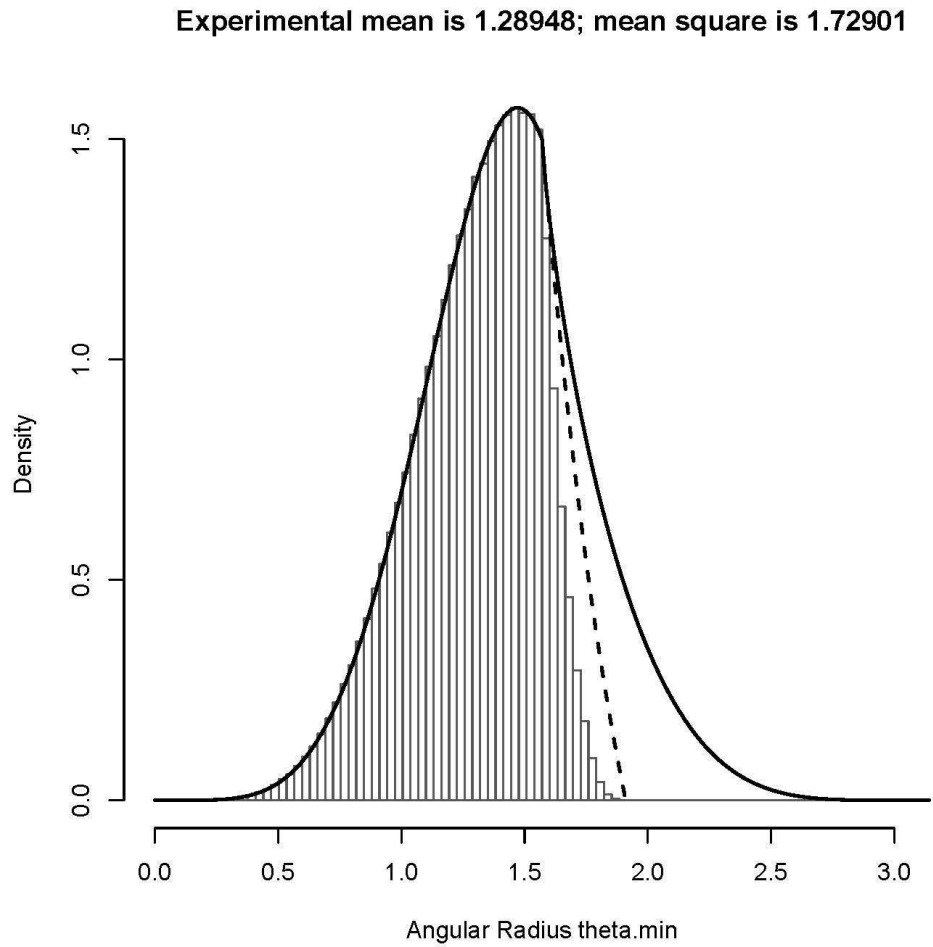}%
\caption{Experimental histogram of $\theta_{\min}$-values and theoretical
fit/upper bounds.}%
\end{figure}

\section{Bounds on Coverage Probability}

Define a function%
\[%
\begin{array}
[c]{ccc}%
q(\omega)=\frac{1}{8}\left(  1-32\kappa\,G(\pi-\omega)\right)  , &  &
0<\omega<\frac{\pi}{2}%
\end{array}
\]
then $q(\omega)$ is a nontrivial lower bound for $p(\omega)$ when
$q(\omega)>0$, that is, when $\omega>84.25^{\circ}$. As discussed earlier,
$q(88^{\circ})=0.0765$. \ 

Our method fails to improve upon the inequality $p(\omega)\geq0$ for
$\omega<84.25^{\circ}$. \ We assess the effect of more closely approximating
the density $\varphi$. \ Let \
\[
\theta_{0}=\pi-\omega_{0}=\arccos\left(  -\tfrac{1}{3}\right)
=1.91063323...\approx109.47^{\circ}.
\]
Assuming the function $1+\varphi$ is logarithmically convex on $[\pi
/2,\theta_{0}]$, it follows that%
\[
\psi_{\operatorname{lcv}}(\theta)=\left(  \frac{5}{2}\right)  ^{\frac
{\theta_{0}-\theta}{\theta_{0}-\pi/2}}-1
\]
dominates $\varphi$. \ This is suggested in Figure 2, where the curve
$\psi_{\operatorname{lcv}}$ is represented by the dashed line. Integrating
$\psi_{\operatorname{lcv}}$ from $\pi/2$ to $\theta$, we obtain $\Psi
_{\operatorname{lcv}}$, which in turn is used to define \
\[%
\begin{array}
[c]{ccc}%
q_{\operatorname{lcv}}(\omega)=\frac{1}{8}-\Psi_{\operatorname{lcv}}%
(\pi-\omega), &  & 0<\omega<\frac{\pi}{2}.
\end{array}
\]
This is a tighter bound than the preceding; for example,
$q_{\operatorname{lcv}}(88^{\circ})=0.0766$. It is nontrivial when
$q_{\operatorname{lcv}}(\omega)>0$, that is, when $\omega>83.90^{\circ}$. Thus
the effect is only slight. \ A rigorous proof that log convexity holds,
however, is not known. \ 

A sharper upper bound for $p(\omega)$ than that provided in \cite{Gbt} is also desired.

\section{Addendum}

The preceding text was written in 2011. \ A question, \textquotedblleft Does
$\kappa$ possess a closed-form expression?\textquotedblright\ appeared in
\cite{WC1}; an affirmative answer
\[
\kappa=\frac{11}{96}-\frac{1}{8\pi^{2}}
\]
appeared in \cite{WC2} with a different line of reasoning. We acknowledge that
the preceding text is heuristic and non-rigorous. \ Sample open issues
include:\ Must a well-centered tetrahedron $ABCD$ possess at least one
acute-angled triangular face?\ \ (Simulation suggests, in fact, that there are
at least two.) \ If all faces $ABC$, $ABD$, $ACD$, $BCD$ are acute, then
\[
\theta_{\min}=\min\left\{  \theta_{\operatorname*{abc}},\theta
_{\operatorname*{abd}},\theta_{\operatorname*{acd}},\theta
_{\operatorname*{bcd}}\right\}
\]
and the inequality $4G(\xi)\geq F(\xi)$ becomes clear. \ It is possible,
however, that only two or three faces are acute. \ Let $N\in\{2,3,4\}$ be the
number of acute faces, given a random well-centered tetrahedron $ABCD$. \ What
precisely is the distribution of $N$? \ All we know is $E(N)=$ $32\kappa$. How
is the inequality $(32\kappa)G(\xi)\geq F(\xi)$ proved? \ (The apparent
continuity of $\psi(\theta)$ at $\theta=\pi/2$ gives only an impression that
this may be correct.) \ Discussion would be appreciated.

\end{document}